\documentclass[11pt]{article}
\usepackage{tipa}
\usepackage[all,import]{xy}
\usepackage{graphicx,color}
\usepackage{amsmath}
\usepackage{amsbsy}
\usepackage{amssymb}
\usepackage{cases}
\usepackage{enumerate}
\usepackage{bm}
\usepackage[colorlinks,linkcolor=black,anchorcolor= black,citecolor= black,urlcolor=black]{hyperref}
\RequirePackage{hyperref}
\RequirePackage{hypernat}

\usepackage{amsthm}
\usepackage{multirow}
\usepackage{subfigure}
\usepackage{colortbl}
\usepackage{algorithm}
\usepackage{algorithmic}

\usepackage{times}
\usepackage{mathptmx}

\definecolor{mygray}{gray}{.9}

\makeatletter
\newcommand*{\rom}[1]{\expandafter\@slowromancap\romannumeral #1@}
\makeatother

\begin{document}

\theoremstyle{definition}
\newtheorem{assumption}{Assumption}
\newtheorem{theorem}{Theorem}
\newtheorem{lemma}{Lemma}
\newtheorem{example}{Example}
\newtheorem{definition}{Definition}
\newtheorem{corollary}{Corollary}

\def\letas{\mathrel{\mathop{=}\limits^{\triangle}}}
\def\ind{\begin{picture}(9,8)
         \put(0,0){\line(1,0){9}}
         \put(3,0){\line(0,1){8}}
         \put(6,0){\line(0,1){8}}
         \end{picture}
        }
\def\nind{\begin{picture}(9,8)
         \put(0,0){\line(1,0){9}}
         \put(3,0){\line(0,1){8}}
         \put(6,0){\line(0,1){8}}
         \put(1,0){{\it /}}
         \end{picture}
    }

\def\AVar{\text{AsyVar}}
\def\Var{\text{Var}}
\def\Cov{\text{Cov}}
\def\sumn{\sum_{i=1}^n}
\def\summ{\sum_{j=1}^m}
\def\convergeas{\stackrel{a.s.}{\longrightarrow}}
\def\converged{\stackrel{d}{\longrightarrow}}
\def\iidsim{\stackrel{i.i.d.}{\sim}}
\def\indsim{\stackrel{ind}{\sim}}
\def\asim{\stackrel{a}{\sim}}
\def\d{\text{d}}

\setlength{\baselineskip}{2\baselineskip}

\setcounter{MaxMatrixCols}{20}

\title{\bf \huge On the Conditional Distribution of the Multivariate $t$ Distribution}
\author{Peng Ding\footnote{
Peng Ding is Assistant Professor (Email: \url{pengdingpku@gmail.com}), Department of Statistics, University of California, Berkeley, 425 Evans Hall, Berkeley, CA 94720, USA. The author thanks the Editor, Associate Editor and three reviewers for their helpful comments. 
}
}

\date{}
\maketitle

\begin{abstract}
As alternatives to the normal distributions, $t$ distributions are widely applied in robust analysis for data with outliers or heavy tails. The properties of the multivariate $t$ distribution are well documented in Kotz and Nadarajah's book, which, however, states a wrong conclusion about the conditional distribution of the multivariate $t$ distribution. Previous literature has recognized that the conditional distribution of the multivariate $t$ distribution also follows the multivariate $t$ distribution. We provide an intuitive proof without directly manipulating the complicated density function of the multivariate $t$ distribution.

\end{abstract}

\section{Introduction}
The conventional version of the multivariate $t$ (MVT) distribution $ \bm{X} \sim  \bm{t}_p(\bm{\mu}, \bm{\Sigma}, \nu)$, with location $\bm{\mu}$, scale matrix $\bm{\Sigma}$, and degrees of freedom $\nu$, has the probability density function
\begin{eqnarray}
\label{eq::pdf}
f(\bm{x}) = 
\frac{ \Gamma \{  (\nu + p)/2  \}  }
{   \Gamma(\nu/2) (\nu \pi) ^{p/2}   |  \bm{\Sigma} |^{1/2}  }
 \left\{  1+ \nu^{-1}   (\bm{x}-\bm{\mu})^\top \bm{\Sigma}^{-1} (\bm{x}-\bm{\mu})    \right\}^{  -(\nu+p)/2 }
.
\end{eqnarray}
We can see from the above that the tail probability of the MVT distribution decays at a polynomial rate, resulting in heavier tails than the multivariate normal distribution. Because of this property, the MVT distribution is widely applied in robust data analysis including linear and nonlinear regressions (Lange et al. 1988; Liu 1994; Liu 2004), linear mixed effects models (Pinheiro et al. 2001), and sample selection models (Marchenko and Genton 2012; Ding 2014).

In the following discussion, we use $W\sim \chi^2_b/c$ to denote the scaled $\chi^2$ distribution, with density proportional to $ w^{b/2-1} e^{-cw/2}.$
We exploit the following representation of the MVT distribution: 
\begin{eqnarray}
\label{eq::def}
\bm{X} = \bm{\mu} + \bm{\Sigma}^{1/2} \bm{Z} / \sqrt{q} ,
\end{eqnarray} 
where $\bm{Z}$ follows a $p$ dimensional standard normal distribution, $  q \sim \chi^2_\nu /  \nu$, and $\bm{Z}$ is independent of $ q$ (Kotz and Nadarajah 2004; Nadarajah and Kotz 2005). It differs from the multivariate normal distribution $\mathcal{N}_p(  \bm{\mu}, \bm{\Sigma} )$ only by the random scaling factor $\sqrt{q}$.

The above representation (\ref{eq::def}) implies $\bm{X}\mid q\sim \mathcal{N}_p(   \bm{\mu} ,  \bm{\Sigma}/q )$, i.e., $\bm{X}$ follows a multivariate normal distribution given the latent variable $q.$
If we partition $\bm{X}$ into two parts, $\bm{X}_1$ and $\bm{X}_2$, with dimensions $p_1$ and $p_2$, we obtain the following normal mixture representation conditional on $q$:
\begin{equation}
\label{eq::rep}
\bm{X} = \begin{pmatrix}
\bm{X}_1\\
\bm{X}_2
\end{pmatrix}  \mid  q
\sim \mathcal{N}_{p_1 + p_2} \left\{  \begin{pmatrix}
\bm{\mu}_1\\
\bm{\mu}_2
\end{pmatrix},
\begin{pmatrix}
\bm{\Sigma}_{11} & \bm{\Sigma}_{12}\\
\bm{\Sigma}_{21} & \bm{\Sigma}_{22}
\end{pmatrix}/q \right\},
\end{equation}
where the location and scale parameters are partitioned corresponding to the partition of $\bm{X}$.
Marginally, we have $\bm{X}_1\mid q \sim \mathcal{N}_{p_1} (  \bm{\mu}_1 ,  \bm{\Sigma}_{11} /q ) $, and therefore
$$
\bm{X}_1   \sim \bm{t}_{p_1}(\bm{\mu}_1, \bm{\Sigma}_{11}, \nu),
$$ 
which follows a $p_1$ dimensional MVT distribution with degrees of freedom $\nu$.

Although we can obtain the marginal distribution in an obvious way, the conditional distribution of $\bm{X}_2$ given $\bm{X}_1$ is less transparent. In fact, the conditional distribution of the MVT distribution is also a MVT distribution, with degrees of freedom different from the original distribution.

\section{Conditional Distribution}
\label{sec::conditional}

Kotz and Nadarajah (2004, page 17) and Nadarajah and Kotz (2005) claimed that the conditional distribution of the MVT distribution is not a MVT distribution except for some extreme cases. For the case with $\bm{\mu}=0$, define $\bm{x}_{2|1} = \bm{x}_1 - \bm{\Sigma}_{21}\bm{\Sigma}_{11}^{-1}\bm{x}_1$, and define $  \bm{\Sigma}_{22|1} =  \bm{\Sigma}_{22} - \bm{\Sigma}_{21}\bm{\Sigma}_{11}^{-1}\bm{\Sigma}_{12}$ as the Schur complement of the block $\bm{\Sigma}_{11}$ in matrix $\bm{\Sigma}.$ They derived the conditional density of $\bm{X}_2$ given $ \bm{X}_1$ by calculating $f(\bm{x}) / f_1(\bm{x}_1)$ using the probability density function in (\ref{eq::pdf}). Ignoring the normalizing constant in formula (15) of Nadarajah and Kotz (2005), we present only the key term:
\begin{eqnarray}\label{eq::conditional-pdf}
f_{2|1}(\bm{x}_2\mid \bm{x}_1)
\propto 
\left\{    
1 +  (\nu + p_1)^{-1}  \bm{x}_{2|1}^\top       
 \left(   \frac{  \nu+ \bm{x}_1^\top \bm{\Sigma}_{11}^{-1} \bm{x}_1  }{ \nu + p_1 }      \bm{\Sigma}_{22|1}    \right)^{-1}      
       \bm{x}_{2|1}
\right\}^{-(\nu+p_1+p_2)/2} .
\end{eqnarray}
Kotz and Nadarajah (2004, page 17) and Nadarajah and Kotz (2005) obtained the correct conditional density function, but due to the complex form of their formula (15), they did not recognize that the key term in (\ref{eq::conditional-pdf}) is essentially the unnormalized probability density function of a $p_2$ dimensional MVT distribution with location $\bm{\Sigma}_{21}\bm{\Sigma}_{11}^{-1}\bm{x}_1$, scale matrix $ ( \nu+ \bm{x}_1^\top \bm{\Sigma}_{11}^{-1} \bm{x}_1  )/( \nu + p_1)\times      \bm{\Sigma}_{22|1}   $, and degrees of freedom $ ( \nu+p_1).$

In the literature, some authors stated the right conclusion about the conditional distribution without providing a proof (e.g., Liu 1994, page 5, DeGroot 2005, page 61), and some other authors re-derived the correct conditional density function in (\ref{eq::conditional-pdf}) and pointed out that it is another MVT distribution (e.g., Roth 2013b). Instead of directly calculating the conditional density or using the general theory of elliptically contoured distributions (Cambanis et al. 1981; Fang et al. 1990, Theorems 2.18 and 3.8; Kibria and Joarder 2006), we provide an elementary and thus  more transparent proof of the conditional distribution based on the normal mixture representations in (\ref{eq::rep}).

%

\section{Proof via Representation}
\label{sec::representation}

Our proof proceeds in two steps: we first condition on $q$, and then average over $q$.

First, we condition on both $\bm{X}_1$ and $q$. According to the property of the multivariate normal distribution, $\bm{X}_2$ follows a multivariate normal distribution conditional on $(\bm{X}_1, q)$, i.e.,
\begin{eqnarray*}
\bm{X}_2\mid( \bm{X}_1, q) 
\sim \mathcal{N}_{p_2}\left\{  
\bm{\mu}_2 + \bm{\Sigma}_{21} \bm{\Sigma}_{11}^{-1} (\bm{X}_1 - \bm{\mu}_1), 
(\bm{\Sigma}_{22} - \bm{\Sigma}_{21}\bm{\Sigma}_{11}^{-1}\bm{\Sigma}_{12})/q
\right\}  
= \mathcal{N}_{p_2}\left(  
\bm{\mu}_{2|1} , \bm{\Sigma}_{22|1} /q
\right),
\end{eqnarray*}
where $\bm{\mu}_{2|1}  =  \bm{\mu}_2 + \bm{\Sigma}_{21} \bm{\Sigma}_{11}^{-1} (\bm{X}_1 - \bm{\mu}_1)$ is the linear regression of $\bm{X}_2$ on $\bm{X}_1$.

Second, we obtain the conditional distribution of $q$ given $\bm{X}_1$. According to Bayes' Theorem, the conditional probability density function of $q$ given $\bm{X}_1$ is 
\begin{eqnarray*}
f_{q|1}(q\mid \bm{X}_1)  &\propto &  f_{1|q}(\bm{X}_1\mid q) \pi(q)\\
 &\propto&  | \bm{\Sigma}_{11} / q  |^{-1/2}   \exp\left\{   -\frac{q}{2} (   \bm{X}_1 - \bm{\mu}_1)^\top \bm{\Sigma}_{11}^{-1}  (   \bm{X}_1 - \bm{\mu}_1)      \right\}
q^{\nu/2 - 1} e^{-\nu q/2} 
\\
&\propto & q^{(\nu + p_1)/2 - 1}   \exp\left\{   -\frac{q}{2}  \left( \nu + d_1 \right)    \right\}, 
\end{eqnarray*}
where $d_1 =    (   \bm{X}_1 - \bm{\mu}_1)^\top \bm{\Sigma}_{11}^{-1}  (   \bm{X}_1 - \bm{\mu}_1)$ is the squared Mahalanobis distance of $\bm{X}_1$ from $\bm{\mu}_1$ with scale matrix $\bm{\Sigma}_{11}$. Therefore, the conditional distribution of $q$ given $\bm{X}_1$ is
$
q\mid \bm{X}_1 \sim  \chi^2_{\nu + p_1}  /  (  \nu +     d_1  ),
$
according to the density of the scaled $\chi^2$ distribution.

Finally, we can represent the conditional distribution of $\bm{X}_2$ given $\bm{X}_1$ as
\begin{eqnarray}
\bm{X}_2\mid \bm{X}_1 
&\sim&
\left(  \bm{\mu}_{2|1}  +  \bm{\Sigma}_{22|1}^{1/2}   \bm{Z}_2 / \sqrt{q} \right) \mid \bm{X}_1 \nonumber \\ 
&\sim & 
\bm{\mu}_{2|1}  +  \bm{\Sigma}_{22|1}^{1/2}   \bm{Z}_2 \Big/ \sqrt{        \frac{  \chi^2_{\nu + p_1} }{ \nu +   d_1     }          } \nonumber  \\
&\sim& 
\bm{\mu}_{2|1} 
+ \left(   \sqrt{      \frac{ \nu +   d_1     }{  \nu + p_1 }        }   \bm{\Sigma}_{22|1}^{1/2} \right) 
\left(   \bm{Z}_2  \Big/ \sqrt{   \frac{  \chi^2_{\nu + p_1}}{ \nu + p_1 }     } \right) ,\label{eq::rep2}
\end{eqnarray}
where $\bm{Z}_2$ is a $p_2$ dimensional standard multivariate normal vector, independent of $ \chi^2_{\nu + p_1} / (\nu+p_1).$
From the normal mixture representation of the MVT distribution in (\ref{eq::def}), we can obtain the following conditional distribution.

\noindent {\bfseries Conclusion One}:
The conditional distribution of $\bm{X}_2$ given $\bm{X}_1$ is
$$
\bm{X}_2\mid \bm{X}_1 \sim \bm{t}_{p_2}  \left(   \bm{\mu}_{2|1} , 
  \frac{ \nu +   d_1     }{  \nu + p_1 } 
  \bm{\Sigma}_{22|1} , 
\nu + p_1
       \right).
$$

The conditional distribution of the multivariate $t$ distribution is very similar to that of the multivariate normal distribution. The conditional location parameter is the linear regression of $\bm{X}_2$ on $\bm{X}_1$. The conditional scale matrix is $\bm{\Sigma}_{22|1}$ inflated or deflated by the factor $(\nu+d_1) / (\nu + p_1)$. 
Because $d_1\sim p_1 \times F(p_1, \nu)$ with mean being $ p_1\nu/(\nu-2) $ or infinity according to $\nu>2$ or $\nu\leq 2$ (Roth 2013a, page 83), on average the factor $(\nu+d_1) / (\nu + p_1)$ is larger than one. With more extreme values of $\bm{X}_1$, the conditional distributions of $\bm{X}_2$ are more disperse.     
More interestingly, the conditional degrees of freedom increase to $(\nu+p_1)$. The more dimensions we condition on, the less heavy-tailedness we have. When $\nu=1$, the MVT distribution is also called the multivariate Cauchy distribution. The marginal distributions are still Cauchy, but the conditional distributions are MVT instead of Cauchy. Although the marginal means of the multivariate Cauchy distribution do not exist for any dimensions, all the conditional means exist because the conditional distributions of  the multivariate Cauchy distribution follow MVT distributions with degrees of freedom at least as large as two.

As a byproduct, we can easily obtain from the representation (\ref{eq::rep2}) that
$$
\sqrt{  \frac{  \nu + p_1 }{ \nu +  d_1   } }
 \left(  \bm{X}_2    -    \bm{\mu}_{2|1}  \right) \mid \bm{X}_1
 \sim \bm{t}_{p_2} \left(    \bm{0},   \bm{\Sigma}_{22|1} , \nu + p_1  \right) ,
$$
which does not depend on $\bm{X}_1$. The above result further implies the following conclusion.

\noindent {\bfseries Conclusion Two}:
 $\bm{X}_1  $ and 
$
\sqrt{  \frac{  \nu + p_1 }{ \nu +  d_1     } }
 \left(  \bm{X}_2    -   \bm{\mu}_{2|1}   \right)
 $
are independent.

It is straightforward to show that $\bm{X}_1  $ and 
$
 \left(  \bm{X}_2    -   \bm{\mu}_{2|1}   \right)
 $
are uncorrelated for both multivariate normal and MVT distributions. If $\bm{X}$ follows the multivariate normal distribution, this also implies independence between $\bm{X}_1  $ and 
$
 \left(  \bm{X}_2    -   \bm{\mu}_{2|1}   \right)
 $. If $\bm{X}$ follows the MVT distribution, we need only to adjust the linear regression residual $\left(  \bm{X}_2    -   \bm{\mu}_{2|1}  \right)$  by the factor $\sqrt{  ( \nu + p_1 ) / ( \nu +  d_1     ) }$ and the independence result still holds.

\section{Discussion}

Kotz and Nadarajah (2004, page 17) and Nadarajah and Kotz (2005) failed to recognize that the conditional distribution of the MVT distribution is also a MVT distribution due to the complexity of the conditional density function. Our proof, based on the normal mixture representation of the MVT distribution, offers a more direct and transparent way to revealing the property of the conditional distribution of the MVT distribution. Our representation in Section \ref{sec::representation} implicitly exploits the data augmentation formula $f_{2|1}(\bm{x}_2\mid \bm{x}_1) = \int f_{2|1,q}(\bm{x}_2\mid \bm{x}_1, q)  f_{q|1}(q\mid \bm{x}_1) dq$ (Tanner and Wong 1987). To make our proof more intuitive, we used formula (5) to avoid the integral, which is essentially a Monte Carlo version of the data augmentation formula.

Conclusion One could also follow from the general theory of elliptically contoured distributions, including the MVT distribution as a special case (Cambanis et al. 1980; Fang et al. 1990). Conditional distributions of elliptically contoured distributions are also elliptically contoured distributions. But this does not immediately guarantee that conditional distributions of the MVT distributions are also MVT distributions without some further algebra. The MVT distribution has the property of having MVT marginal and conditional distributions. It will be interesting to find other sub-classes of elliptically contoured distributions that have the same property.


\section*{References}
\begin{description} \itemsep=-\parsep \itemindent=-1.2 cm

\item
Cambanis, S., Huang, S. and Simons, G. (1981). On the theory of elliptically contoured distributions. {\it Journal of Multivariate Analysis}, {\bfseries 11}, 368--385.

\item 
Ding, P. (2014). Bayesian robust inference of sample selection using selection-$t$ models. {\it Journal of Multivariate Analysis}, {\bfseries 124}, 451--464.

\item
DeGroot, M. H. (2005). {\it Optimal Statistical Decisions}. New York: John Wiley \& Sons.

\item
Fang, K. T., Kotz, S., and Ng, K. W. (1990). {\it Symmetric Multivariate and Related Distributions}. London: Chapman \& Hall.

\item
Gelman, A., Carlin, J. B., Stern, H. S., Dunson, D., Vehtari, A., and Rubin, D. B. (2014). {\it Bayesian Data Analysis, 2nd edition}. London: Chapman \& Hall/CRC.

\item
Kibria, B. G., and Joarder, A. H. (2006). A short review of multivariate $t$-distribution. {\it Journal of Statistical Research}, {\bfseries 40}, 59--72.

\item
Kotz, S. and Nadarajah, S. (2004). {\it Multivariate $t$ Distributions and Their Applications}. New York: Cambridge University Press.

\item
Lange, K. L., Little, R. J., and Taylor, J. M. (1989). Robust statistical modeling using the $t$ distribution. {\it Journal of the American Statistical Association}, {\bfseries 84}, 881--896.

\item
Liu, C. (1994). {\it Statistical Analysis Using the Multivariate t Distribution}. Doctoral dissertation, Harvard University. 

\item
Liu, C. (2004). Robit regression: a simple robust alternative to logistic and probit regression. In Gelman, A. and Meng, X.-L. (Eds.) {\it Applied Bayesian Modeling and Causal Inference from Incomplete-Data Perspectives} (pp. 227--238). New York: John Wiley \& Sons.

\item 
Marchenko, Y. V., and Genton, M. G. (2012). A Heckman selection-$t$ model. {\it Journal of the American Statistical Association}, {\bfseries 107}, 304--317.

\item 
Nadarajah, S. and Kotz, S. (2005). Mathematical properties of the multivariate $t$ distribution. {\it Acta Applicandae Mathematica}, {\bfseries 89}, 53--84.

\item
Pinheiro, J. C., Liu, C., and Wu, Y. N. (2001). Efficient algorithms for robust estimation in linear mixed-effects models using the multivariate $t$ distribution. {\it Journal of Computational and Graphical Statistics}, {\bfseries 10}, 249--276.

\item
Roth, M. (2013a). {\it Kalman Filters for Nonlinear Systems and Heavy-Tailed Noise}, Licentiate of  Engineering Thesis, Link\"oping University, Link\"oping.

\item
Roth, M. (2013b). On the multivariate $t$ distribution. Technical report from Automatic Control at Link\"opings Universitet. Report no.: LiTH-ISY-R-3059. 

\item
Tanner, M. A., and Wong, W. H. (1987). The calculation of posterior distributions by data augmentation. {\it Journal of the American Statistical Association}, {\bfseries 82}, 528--540.

\end{description}

\end{document}